\documentclass[12pt,twoside,reqno]{amsart}
\usepackage[colorlinks=true,citecolor=blue]{hyperref}
\usepackage{mathptmx, amsmath, amssymb, amsfonts, amsthm, mathptmx, enumerate, color,mathrsfs,paralist}
\usepackage{graphicx}
 \usepackage{color}
\usepackage{tikz}

\usepackage{epstopdf}
\usepackage[pagewise]{lineno}

\newtheorem{theorem}{Theorem}[section]

\theoremstyle{definition}

\newtheorem{remark}[theorem]{Remark}
\numberwithin{equation}{section}
\newcommand{\be}{\begin{equation}}
\newcommand{\ee}{\end{equation}}
\newcommand{\la}{\lambda}
\newcommand{\lap}{\mbox{$\Delta$}}
\allowdisplaybreaks[4]

\begin{document}
\setcounter{page}{1}

\vspace*{2.0cm}
\title[Symmetry and monotonicity of solutions ]
{Symmetry and monotonicity of solutions to elliptic and   parabolic fractional $p$-equations }
\author[P. Wang]{Pengyan Wang$^{*}$  }
\maketitle
\vspace*{-0.6cm}

\begin{center}
{\footnotesize

 School of   Mathematics and Statistics,
Xinyang Normal University,
  Xinyang, Henan,   464000, China

}\end{center}

\vskip 4mm {\footnotesize \noindent {\bf Abstract.}
In this article we   first establish the maximum principle of the antisymmetric functions for  parabolic fractional $p$-equations. Then we use it and the parabolic inequalities to provide a  different proof of symmetry and monotonicity for solutions to  elliptic  fractional $p$-equations with gradient terms.
 Finally, base on suitable initial value, by the maximum principle of the antisymmetric functions for  parabolic fractional $p$-equations,  we  attain symmetry and monotonicity of positive solutions in each finite time to  nonlinear parabolic fractional $p$-equations  on the whole space and   bounded domains.

 We believe that the maximum principle and parabolic inequalities
obtained here   can be utilized to many elliptic and parabolic problems    involving nonlinear nonlocal operators.

 \noindent {\bf Keywords.}
The parabolic fractional $p$-equations, monotonicity,  symmetry, the method of moving planes,   elliptic  fractional $p$-equations.

 \noindent {\bf 2010 Mathematics Subject Classification.}
35K58,  35A09, 35B06, 35B09.}

\renewcommand{\thefootnote}{}
\footnotetext{ $^*$Corresponding author.
\par
E-mail addresses: wangpy@xynu.edu.cn (P. Wang).
\par
  }

\section{Introduction}

Symmetry and monotonicity of positive solutions of the local elliptic equations in unit ball was first obtained
by Gidas, Ni and Nirenberg \cite{GNN}. In recent decades, the elliptic equation  involving nonlocal  operators, notably fractional Laplace operators,  have received extensive attention and have achieved a large number of results (see \cite{CL,CLL,CLZ,CZ1,DL,DLi,L1,LZ2})   since
the work of Caffarelli and Silvestre \cite{CS}. For other results on fractional Laplacian and fractional $p$-Laplacian
equations, we refer readers to \cite{CT} for regularity  and the maximum principle,    \cite{DPV,WCN} for existence
and symmetry results of the Schr\"{o}dinger  equation, \cite{OS} for regularity up to the
boundary,  \cite{BPSV} for a widder's type Theorem,
  \cite{SV} for mountain pass solutions, \cite{L,WC2} for monotonicity of  solutions, \cite{CLM,CQ} for a review, and references therein.

In this paper,  we first  prove  a  maximum principle   for antisymmetric functions to parabolic fractional $p$-equations. Then we
    investigate the following elliptic fractional $p$-equation with a gradient term
\begin{equation} \label{eq:n20191}
\left\{\begin{array}{lll}
(-\lap)^s_p u(x)= h(x,u(x), \nabla u(x)),& x\in \Omega,\\
u(x)>0,& x\in \Omega,\\
u(x)=0,& x\in  \Omega ^c,
\end{array}\right.
\end{equation}
where $\Omega$ is a bounded domain in $\mathbb R^n$ which is convex in $x_1$ direction and the fractional $p$-Laplacian with $ s\in (0,1),~  p\in[2,+\infty)$  is given by $$
(-\lap)^s _p u(x ) =C_{n,sp} P.V. \int_{\mathbb R^n} \frac{ |u(x )-u(y )|^{p-2}(u(x )-u(y )) }{|x-y|^{n+sp}}dy,
$$
where $P.V.$ stands for the Cauchy principal value. When $p=2$, $(-\lap)^s_p$ becomes the usual fractional Laplacian $(-\lap)^s$.
When $s=1$, $(-\lap)^s_p$ becomes the well-known $p$-Laplacian $-\lap_p$. When $s=1,p=2$, $(-\lap)^s_p$  becomes standard Laplacian $-\lap$.
Define $$   L _{sp}=\{u  \in L^{p-1}_{loc} (\mathbb{R}^n) \mid \int_{\mathbb R^n} \frac{|1+u(x )|^{p-1}}{1+|x|^{n+sp}}dx<+\infty\},$$
then it is easy to see that for $u\in C^{1,1}_{loc} \cap   L _{sp}$, $(-\lap)^s_p u(x ) $ is well defined.
In the following, we denote $\nabla u=(\frac{\partial u}{\partial x_1}  ,\frac{\partial u}{\partial x_2} , \cdots, \frac{\partial u}{\partial x_n} ) $ by $\textbf{q}= (q_1,q_2,\cdots,q_n)$ and prove
\begin{theorem}\label{eqn22}
Let $u(x)\in C^{1,1}_{loc}(\Omega) \cap C(\bar \Omega)$ be a solution of   elliptic fractional  $p$-equation \eqref{eq:n20191}, where $g(x,u,\textbf{q})$ is Lipschitz continuous with respect to $(u,\textbf{q})$ and
\begin{equation}\label{eq:n288}\aligned
&\mbox{for}~x_1<\tau_1,~\tau_1+x_1<0~\mbox{and}~q_1\geq 0,\\
& h(x_1,x',u,q_1,q_2,\cdots,q_n)\leq h(\tau_1,x',u,-q_1,q_2,\cdots,q_n).
\endaligned\end{equation}
Then $u(x_1, x')$ is strictly increasing in the left half of $\Omega$
 in $x_1$-direction and
\begin{equation}\label{eq:n299}
u(x_1,x')\leq u(\tau_1,x'),~\forall~x_1<\tau_1,~x_1+\tau_1<0,~(x_1,x') \in \Omega.
\end{equation}
Furthermore if $h(x_1,x',u,q_1,q_2,\cdots,q_n)=h(-x_1,x',u,-q_1,q_2,\cdots,q_n)$,  then
$$ u(x_1,x')=u(-x_1,x').$$
\end{theorem}
\begin{remark}
 Theorem \ref{eqn22}  has been obtained by the method of moving planes and the maximum principle for fractional elliptic equations in \cite{CHL} with $p=2$ and \cite{DLW} with $p\geq2$. The monotonicity along the $x_1$ direction in the entire domain $\Omega$ has been obtained in \cite{WL} by a sliding method with $p=2$, and in \cite{WCD,W} with $p\geq 2$. Here we apply the parabolic inequality and  a  maximum principle of antisymmetric functions for parabolic fractional $p$-equations  to give   a completely  different proof. Our result is to extend   the result in \cite{Z} to the fractional $p$-equation.  These results can be seen as the nonlinear nonlocal counterpart  to some of the results in the celebrated article \cite{BN1}. This indicates that the elliptic equation is closely related to its corresponding parabolic equation.
\end{remark}
   Due to the nonlinearity of fractional  $p$-Laplacian, its research is more difficult than that of fractional  Laplacian. Secondly, the presence of gradient terms in the nonlinear term at the right end of the equation \eqref{eq:n20191} adds difficulty to our research. To overcome these difficulties,   we first establish a maximum principle for fractional $p$-Laplacian and derive parabolic inequalities using nonlinear functions with gradient terms and Lipschitz continuity. Then, in the proof of using the moving plane method, as in \cite{CHL}, when deriving the contradiction for the minimum point of the function $w_\la(x)$(as defined in Section 2 below), we  apply the approach of simultaneously finding the minimum value of   $w_\la(x)$ for  both variables $\la$ and $x$.
   Therefore, this method is effective for equations with nonlinear terms containing gradient terms.

 There are also many results regarding the symmetry and monotonicity of solutions to parabolic equations, and we only list some relevant literature.   After Dancer and Hess \cite{DH} considered to time periodic
solutions, symmetry of general positive solutions of parabolic equations
on bounded domains was studied in \cite{B,HP,P1,P2}.  In 1989, Li \cite{L2} based on symmetric initial values derived symmetry for positive solutions to fully nonlinear  parabolic equations. Pol\'{a}\v{c}ik \cite{P3} investigated the symmetry  for positive solutions to local quasilinear parabolic equations in the whole space and show that  its  solutions are symmetric at all time $t<T$.  Up to now, there are few results on the symmetry and monotonicity for solutions to parabolic equations involving nonlocal operators.
  In 2012, Jin and Xiong \cite{JX} obtained  monotonicity   for the solution to a fractional Yamabe flow equation  by  a strong maximum principle and a Hopf lemma. In 2014, Jarohs and Weth \cite{JW1} arrived at the asymptotic symmetry for weak solutions to a nonlinear fractional reaction-diffusion equation  in the bounded domain. Very recently, Chen  et al. \cite{CWNH} developed an asymptotic method of moving planes and derived asymptotic symmetry for solutions to fractional  parabolic  equations. Zeng \cite{Z} obtained symmetry of solutions to the fractional  parabolic  equation  with symmetric initial value for each finite time.

Another purpose of this paper is to investigate the symmetry and monotonicity of the solution  to the following nonlinear parabolic fractional $p$-equation
\begin{equation}\label{eq:n1}
\left\{\begin{array}{lll}
\frac{\partial u}{\partial t}+(-\lap)^s_pu(x,t)=F(x,t,u(x,t)),& (x,t)\in \Omega \times (0,T],\\
u(x,t)=\Phi (x,t),& (x,t)\in (\mathbb R^n\backslash \Omega)\times [0,T],
\end{array}\right.
\end{equation}
where    $\Omega \subset \mathbb R^n$ is a bounded  domain  which is symmetric and convex in the $x_1$ direction and $ n\geq 2,~ T >0 $. We call the domain $\Omega$ is convex in $x_1$ direction if and only if $ (x_1,x'),(\tau_1,x')\in \Omega$ imply that
$(\alpha x_1+(1-\alpha)\tau_1,x') \in \Omega$ for $0<\alpha<1$.
For every fixed positive $t $,
 $$
(-\lap)^s _p u(x,t) =C_{n,sp} P.V. \int_{\mathbb R^n} \frac{ |u(x,t)-u(y,t)|^{p-2}(u(x,t)-u(y,t)) }{|x-y|^{n+sp}}dy,
$$
where $  p\in[2,+\infty), s\in (0,1)$ and $P.V.$ stands for the Cauchy principal value.
It is easy to see that    for $u \in C^{1,1}_{loc} (\mathbb R^n)\cap   {\mathcal L}_{sp } ,$ $(-\lap)_p^su$ is well defined, where $$ {\mathcal L}_{sp}=\{u(\cdot, t) \in L^{p-1}_{loc} (\mathbb{R}^n) \mid \int_{\mathbb R^n} \frac{|u(x,t)|^{p-1}}{1+|x|^{n+sp}}dx<+\infty\}.$$

Our  other results are:
\begin{theorem}\label{thmn1} Let $\Omega\subset\mathbb R^n$ be a bounded  domain. Assume that $u(x,t)\in\big(C^{1,1}_{loc}(\Omega)\cap {\mathcal L}_{sp} \cap   C(\mathbb R^n)\big)\times C^1([0,T])$ is a solution of \eqref{eq:n1} which is continuous on $\bar \Omega \times [0,T] $,  $F$ is Lipschitz continuous with respect to $u$ and satisfies
\begin{equation}\label{20191}\aligned
&F(x_1,x',t,u(x,t))\leq F(\hat x_1,x',t,u(x,t)),\\
& \mbox{for}~x_1\leq \hat x_1 \leq -x_1,~(x_1,x')\in \Omega,~t \in (0,T].
\endaligned\end{equation}
Suppose the initial-exterior conditions of $u$ satisfy the following:
\begin{equation}\label{eq:n2}\aligned
&\mbox{the~ initial~ value}~ u(x,0) ~\mbox{satisfies}\\
 & u (x_1,x',0)<u (\hat x_1,x',0)~\mbox{for}~x_1<\hat x_1 <-x_1,~(x_1,x')\in \Omega
\endaligned \end{equation}
and
\begin{equation}\label{eq:n3}\aligned
& \Phi(y_1,y',t)\leq \Phi(z_1,z',t),~y_1<z_1,~y=(y_1,y'),(z_1,z')\in \Omega^c,~z_1<0,\\
& \Phi(y_1,y',t)\leq  u( x_1,x',t),~
  y \in    \Omega^c ,~ x \in \Omega,~ x_1<0,~ t \in [0,T] .
\endaligned \end{equation}
Then $u$ is monotone increasing in $x_1$ for $x_1<0$ and $u(x_1,x',t)\leq u(-x_1,x',t)$ for
$(x_1,x')\in \Omega,~x_1<0$ and $  t \in [0,  T]$.

In addition, assume that
\begin{equation}\label{w20191}
F(x_1,x',t,u(x,t))= F(- x_1,x',t,u(x,t)),~~\mbox{for}~(x_1,x')\in \Omega, ~t \in (0,T],
\end{equation}
and the  initial-exterior conditions of $u$ are symmetric
in $\{x_1=0\}$, then $u$ is symmetric in $ x_1 $ and has only one crest. That is, for $(x_1,x')\in \Omega,~x_1<0,~ t\in (0,T]$, we have
\be \label{2q2q}
u(x_1,x',t)=u(-x_1,x',t),  \frac{\partial u}{\partial x_1}(x_1,x',t)\geq 0.
\ee
Furthermore,  we have
\begin{equation}\label{eq:n201912}
\frac{\partial u}{\partial x_1}(x_1,x',t)>0~\mbox{for}~(x_1,x')\in \Omega,~x_1<0,~t \in [0,  T].
\end{equation}
\end{theorem}
\begin{remark}
 When the solution of the fractional parabolic equation is negative, the above conclusion still holds. When $p=2$, our results reduce Theorem 1.1 and Theorem 1.3 in \cite{Z}.
\end{remark}

 When $\Omega$ is a  unit ball,    we obtain
\begin{theorem}\label{thmn3}
 Let $u(x,t)\in \big(C^{1,1}_{loc}(B_1(0)) \cap C(  \overline{B_1(0)})\big) \times C^1([0,T)) $ be a positive solution to
 \begin{equation}\label{1eq:n1}
\frac{\partial u}{\partial t}+(-\lap)_p^su(x,t)=F(x,t,u(x,t)),~ (x,t)\in B_1(0) \times (0,T).
\end{equation}
Assume that $F(x,t,u(x,t))=F(|x|,t,u(x,t))$ is monotone    decreasing in $|x|$ and Lipschitz continuous with respect to $u$ for $t\geq t_0$ and some  time $t_0\geq 0$,  the  exterior condition of $u$   satisfies
 \begin{equation}\label{eq:nll13}
 u(x ,t )=0,~(x,t)\in B^c_1(0)\times [t_0,T).
\end{equation}
   Suppose that
\begin{equation}\label{eq:ll112}
   u  (x, t_0) = u (|x|,t_0) ~\mbox{and}~ u(x, t_0)~\mbox{ is    decreasing  in }~ |x| .
\end{equation}
Then for every time  $t\in(t_0,T),$ $u$ is radially symmetry and monotone decreasing about origin, $i.e.,$
\begin{equation}\label{eq:n201918}
u(x ,t)=u(|x|,t),~x\in B_1(0),~t_0 < t <T .
\end{equation}
\end{theorem}
\begin{remark}
 From   theorem \ref{thmn3}, it can be inferred that if the solution for \eqref{eq:n1}  at any fixed time $t_0\geq 0$ is radially symmetric, then  the radial symmetry of the solution can be obtained  for all times between   $   t_0$ and $T$.

\end{remark}

When $\Omega$ is  $\mathbb R^n$, we prove
\begin{theorem}\label{thmn4}
Suppose that $u(x,t)\in (C^{1,1}_{loc} (\mathbb R^n)\cap {\mathcal L}_{sp}) \times C^1([0,T])$ is a positive solution of
\begin{equation} \label{eq:wn1}
\frac{\partial u}{\partial t}+(-\lap)^s_p u(x,t)=F(x,t,u(x,t)),~(x,t)\in \mathbb R^n \times (0,T].
\end{equation}
The function $F$ is Lipschitz continuous with respect to $u$   and satisfies
\begin{equation}\label{220191}
F(x_1,x',t,u(x,t))\leq F(\tilde x_1,x',t,u(x,t)),~~\mbox{for}~x_1\leq \tilde x_1 ,~x_1\leq 0,~(x_1,x')\in \mathbb R^n.
\end{equation}
Assume that
\be\label{4283}
\underset{ |x|\rightarrow \infty}{\lim}u(x,t)= 0,~ \mbox{for all }~t\in[0,T]
\ee
and
\begin{equation}\label{xeq:n2}
   u  (x_1,x',  0) < u (y,x', 0) ~\mbox{for}~ x_1\leq y \leq -x_1,x_1< 0~\mbox{ and }~ x' =(x_2,\cdots,x_n) .
\end{equation}
 Then $u$
is monotone increasing in $x_1$ for $x_1<0$ and
$$u(x_1,x',t)\leq u(-x_1,x',t)~\mbox{for}~x_1<0,~t\in [0,T].$$

Furthermore, if $F(x,t,u)=F(|x|,t,u)$, and
 \begin{equation*}\label{xeq:n2q}
   u  (x,  0) = u (|x|, 0) ~\mbox{and}~ u(x,  0)~\mbox{ is  monotone   decreasing  with respect to }~ |x|  ,
\end{equation*}
then for each time $0\leq t \leq  T,~x\in  \mathbb R^n$,
 $$u(x ,t)=u(|x|,t).$$
\end{theorem}

 The  paper is organized as follows. In Section 2,  a maximum principle for   anti-symmetric functions is established.
  In Section 3,  combining the maximum principle, parabolic inequalities with the method of moving planes, we   show that Theorem \ref{eqn22}. Section 4, we prove the  monotonicity and symmetry for solutions to nonlinear   parabolic fractional $p$-equations in the whole space and bounded domain.
\section{Maximum principles for antisymmetry functions}
For simplicity,  choose any direction to be the $x_1$ direction.      Let $T_\lambda=\{x\in \mathbb R^n| x_1=\lambda,~ \lambda \in \mathbb R\} $ be the moving planes,
$\Sigma_\lambda=\{x \in \mathbb R^n\mid x_1<\lambda\}$ be the region to the left of the plane,
$x^\lambda=(2\lambda-x_1,x')$ be the reflection of $x$ about the plane $T_\lambda$. Set
$$ u _\lambda(x,t)=u (x^\lambda,t),  ~w_\lambda (x,t)=u _\lambda(x ,t)
-u ( x,t) . $$
We call a function $w_\lambda(x,t)$ is   antisymmetric function $i.e.$, for any $t $, $w_\lambda(x ,t)$ is  antisymmetric function in $\mathbb R^n$ if and only if
\begin{equation}\label{eq:j2}
w_\lambda(x_1,x_2,\cdots,x_n,t)=-w_\lambda(2\lambda-x_1,x_2,\cdots,x_n,t). \end{equation}

We first give a    simple maximum principle.
\begin{theorem}\label{thmn10}(A simple maximum principle\cite{WpC})
Let $p\geq 2,~\Omega$ be a bounded domain in $\mathbb R^n$.
Assume that
 $$u (x,t),  u_\la (x,t) \in (C^{1,1}_{loc} \cap {\mathcal L}_{sp} )\times  C^{1} (   [0,\infty))$$
  and $w_\la(x,t) $ is lower semi-continuous   about $x$ on  $\bar \Omega$.    If
\begin{equation}\label{1j1000}
\left\{\begin{array}{ll}
\frac{\partial w_\la}{\partial t}  (x,t)+(-\lap)_p^{s}u_\la (x,t) -(-\lap)_p^{s}{u} (x,t) \geq 0,    ~ (x,t)\in \Omega\times ( 0,\infty),\\
w_\la(x,t) \geq 0,  ~ \hspace{4.2cm}~(x,t)\in (\mathbb R^n \backslash \Omega )\times [0,T]     ,\\
w_\la(x,0)\geq 0,   ~ \hspace{6.8cm}~ x\in \Omega,
\end{array}
\right.
\end{equation}
then
\be \label{5201}
w_\la(x,t) \geq 0  ~\mbox{in}~  \bar \Omega \times [0,T],~\forall ~T>0 .
\ee
  Furthermore,    \eqref{5201} hold for unbounded region $\Omega$ if we further assume that for all $t\in [0,T]$
\be \label{52110}
\underset{|x| \rightarrow \infty}{\underline{\lim }} w_\la(x,t)\geq 0.
\ee
  Under the conclusion \eqref{5201}, if $w_\la=0$ at some point $(x_0,t_0)\in \Omega \times (0,T]$, then
\be \label{5211}
w_\la(x,t_0)=0,~a.e.~\mbox{ for}~x\in \mathbb R^n.
\ee
\end{theorem}

Next we   develop   the   maximum principle for the anti-symmetric functions to  parabolic fractional $p$-equations.
\begin{theorem} \label{thmn12} (A maximum principle for   anti-symmetric functions)
Let $p\geq 2,~\Omega$ be a   domain in $\Sigma_\la$.
Assume that $w_\lambda(x,t)\in  (C^{1,1}_{loc}(\Omega) \cap {\mathcal L}_{sp})\times C^{1}( [0,\infty))$  is lower semi-continuous with respect to $x$ on $ \bar{\Omega} $     and satisfies
 \begin{equation}\label{1eq:j3}
\left\{\begin{array}{ll}
\frac{\partial w_\lambda}{ \partial  t} +(-\lap)^{s}_pu_\lambda(x,t)-(-\lap)^{s}_pu (x,t)\geq c(x,t)w_\lambda(x,t),  (x,t)\in   \Omega \times (0,\infty),\\
w_\lambda( x^\lambda, t ) =-w_\lambda(x,t),  ~~~~~~~~~~~~~~~~~~~~~~~~~~~~~~~~~~~~~~~ (x,t) \in   \Sigma_\lambda \times [0,\infty),\\
w_\lambda(x,t) \geq0, ~~~~~~~~~~~~~~~~~~~~~~~~~~~~~~~~~~~~~~~~~~~~~  (x,t) \in    (\Sigma_\lambda \backslash \Omega) \times [0,\infty)  ,\\
 w_\lambda (x,0)\geq 0, ~~~~~~~~~~~~~~~~~~~~~~~~~~~~~~~~~~~~~~~~~~~~~~~~~~~~~~~~~~~~~~ ~~~~ x \in \Omega,
\end{array}
\right.
\end{equation}
where $c(x,t)$ is bounded   from above. \\
  ${\rm (i)}$  If
 $\Omega$ is a bounded  domain,       then we have
\begin{equation}\label{1eq:j4}
w_\lambda(x,t)\geq0 \mbox{ in } \Omega\times [0,T];
\end{equation}
${\rm (ii)}$  if $\Omega$ is unbounded, then the conclusion \eqref{1eq:j4} still holds under the additional condition:
\begin{equation}\label{1eq:jj4}
\underset{|x| \rightarrow \infty}{\underline{\lim}}  w_\lambda(x, t) \geq 0 ~\mbox{uniformly for } t\in [0,T] ;
\end{equation}
${\rm (iii)}$  furthermore,     under the conclusion   \eqref{1eq:j4},  if  $w_\lambda $   attains 0 at some point    $ (x_0,t_0)\in \Omega \times (0,T]$, then
\begin{equation}\label{eq:j6}
w_\lambda(x,t_0) = 0,~a.e.~ \mbox{for}~ x \in \mathbb R^n .
\end{equation}
\end{theorem}
\begin{remark}
 The equivalent expression of  ${\rm (iii)}$ is that under the conclusion   \eqref{1eq:j4}, then for any fixed $t\in (0,T]$, we have
 $ w_\lambda(x,t ) >0, x \in \Omega$ or $w_\lambda(x,t ) =0, a.e. $ for $ x \in \mathbb R^n$.
\end{remark}
\begin{remark}
From the proof process, it can be seen that if $\Omega \subset \Sigma_\la$ is a narrow region or $c(x,t)=0$, Theorem  \ref{thmn12} still holds.
When $p=2$,   the result of parabolic fractional Laplacian has been obtained in \cite{CWNH,Z}.
Here, we only  prove the case of $p>2$, the main difference between  parabolic fractional  Laplacian and parabolic fractional $p$-Laplacian is that the latte is a nonlinear operator, $i.e.,$
$$ (-\lap)_p^{s}u _\la(x,t) -(-\lap)_p^{s} {u} (x,t)\not = (-\lap)_p^{s}w_\la (x,t).$$
\end{remark}
\begin{proof}
${\rm (i)}$ Set $m>0$ as a fixed   constant to be  selected later and
$$ \tilde u_\la(x,t)= e^{-mt} u_\la(x,t),~ \tilde u (x,t)= e^{-mt} u (x,t), ~
  \tilde w (x,t)= e^{-mt} w_\la(x,t).$$
Then from \eqref{1eq:j3} we have
\begin{equation}\label{1} \aligned
&\frac{\partial \tilde w}{\partial t}(x,t) +(-\lap)^s_p \tilde u_\la(x,t)-(-\lap)^s_p \tilde u(x,t) \\
=&-m e^{-mt} w_\la(x,t) +e^{-mt} \frac{ \partial w_\la}{\partial t}(x,t) +e^{-m(p-1)t}[(-\lap)^s_p u_\la (x,t)-(-\lap)^s_p u (x,t)].
\endaligned \end{equation}
We   will obtain  \eqref{1eq:j4}  by proving
\begin{equation}\label{11eq:j4}
\tilde w (x,t)\geq 0,~(x,t) \in  \Omega\times [0,T].
\end{equation}
If  \eqref{11eq:j4}  is not valid,
then  the lower semi-continuity of $w_\lambda(x,t)$ with respect to $x$ on $\bar{\Omega} $ implies that there exists a point
$( \bar x,\bar t )\in\bar{\Omega}\times (0,T]$  such that
$$\tilde w(\bar x ,\bar t) :=  \underset{\bar{\Omega}\times(0,T] }{\min}  \,\tilde w(x,t)<0.$$
Due to the same sign for $\tilde w$ and $  w_\la$,   the third and fourth inequalities  for (\ref{1eq:j3})  indicates that the point  $(\bar x , \bar t)$ is in the interior of $\Omega\times (0,T]$.  Hence
$$\frac{\partial{ \tilde w}} {\partial t}(\bar x ,\bar t)\leq 0 .$$

For simplicity's sake, we denote     $H(z)=|z|^{p-2}z$.  Without doubt, $H(z)$ is a strictly increasing function, and $H'(z)=(p-1)|z|^{p-2}\geq 0$.
According to  the definition of parabolic fractional $p$-Laplacian,   it can be seen that
 \begin{equation}\aligned\label{eq:o100}
&(-\lap)_p^s \tilde u_\la(\bar x ,\bar t)-(-\lap)^s_p \tilde  {u}(\bar x, \bar t) \\
=& C_{n,sp} P.V. \int_{\mathbb R^n} \frac{ H(\tilde u_\la(\bar x,\bar t)-\tilde u_\la( y,\bar t))-H(\tilde  {u}(\bar x,\bar t)- \tilde {u}( y,\bar t))}{|\bar x-y|^{n+sp}}dy \\
=&C_{n,sp} P.V. \int_{\Sigma_\la} \frac{ H(\tilde u_\la(\bar x,\bar t)-\tilde u_\la( y,\bar t))-H(\tilde  {u}(\bar x,\bar t)- \tilde {u}( y,\bar t))}{|\bar x-y|^{n+sp}}dy \\
&+C_{n,sp}\int_{\Sigma_\la} \frac{ H(\tilde u_\la(\bar x,\bar t)-\tilde u( y,\bar t))-H( \tilde {u}(\bar x,\bar t)-\tilde  {u}_\la( y,\bar t))}{|\bar x-y^\la|^{n+sp}}dy \\
=& C_{n,sp} P.V.\int_{\Sigma_\la}[\frac{1}{|\bar x-y|^{n+sp}}-\frac{1}{|\bar x-y^\la|^{n+sp}}]\\
 &\hspace{2cm} \cdot [H(\tilde u_\la(\bar x,\bar t)-\tilde u_\la( y,\bar t))-H( \tilde {u}(\bar x,\bar t)- \tilde {u}( y,\bar t))]dy \\
&+C_{n,sp} \int_{\Sigma_\la} [\frac{H(\tilde u_\la(\bar x,\bar t)-\tilde u_\la(y,\bar t))-H(\tilde u (\bar x,\bar t)-\tilde u ( y,\bar t))}{|\bar x-y^\la|^{n+sp}}\\
&\hspace{1cm} +\frac{H(\tilde u_\la(\bar x,\bar t)-\tilde u (y,\bar t))-H(\tilde u (\bar x,\bar t)-\tilde u_\la(y,\bar t))}{|\bar x-y^\la|^{n+sp}}] dy\\
:=&C_{n,sp}  \{I_1+I_2\}.
\endaligned \end{equation}
For $I_1$, we first notice that
$$\frac{1}{|\bar x-y|^{n+sp}}-\frac{1}{|\bar x-y^\la|^{n+sp}}>0,~\forall~\bar x,y \in \Sigma_\la.$$
While for the second part in the integral,
 we derive
$$H(\tilde u_\la(\bar x, \bar t)-\tilde u_\la(y,\bar t))-H( \tilde {u}(\bar x,\bar t)-\tilde  {u}(y, \bar t))\leq 0~\mbox{but} \not\equiv 0,$$
due to the monotonicity of $H$, $(\bar x, \bar t)$ is a minimum point of $\tilde w_\la$ and the fact that
$$[\tilde u_\la(\bar x, \bar t)-\tilde u_\la(y,\bar t)]-[\tilde  {u}(\bar x,\bar t)- \tilde {u}(y,\bar t)]=\tilde w(\bar x,\bar t)-\tilde w(y, \bar t)\leq 0~\mbox{but} ~\not \equiv 0.$$
It yields that
\begin{equation}\label{eq:o102}
I_1< 0.
\end{equation}
Next to estimate $I_2$, we regroup the terms and apply the mean value theorem to derive
 \begin{equation}\label{eq:o101}\aligned
 I_2=& C_{n,sp} \int_{\Sigma_\la} [\frac{H(\tilde u_\la(\bar x,\bar t)-\tilde u_\la(y,\bar t))-H(\tilde u (\bar x,\bar t)-\tilde u ( y,\bar t)) }{|\bar x-y^\la|^{n+sp}} \\
 &~+\frac{H(\tilde u_\la(\bar x,\bar t)-\tilde u (y,\bar t))-H(\tilde u (\bar x,\bar t)-\tilde u_\la(y,\bar t))}{|\bar x-y^\la|^{n+sp}}]dy\\
 =&C_{n,sp} \int_{\Sigma_\la} [\frac{H(\tilde u_\la(\bar x,\bar t)-\tilde u_\la(y,\bar t))-H(\tilde u (\bar x,\bar t)-\tilde u_\la(y,\bar t)) }{|\bar x-y^\la|^{n+sp}} \\
 &~+\frac{H(\tilde u_\la(\bar x,\bar t)-\tilde u (y,\bar t))-H(\tilde u (\bar x,\bar t)-\tilde u ( y,\bar t))}{|\bar x-y^\la|^{n+sp}}]dy\\
=& C_{n,sp}\tilde w_\la(\bar x,\bar t) \int_{\Sigma_\la} \frac{H'(\xi(y,\bar t))+H'(\eta(y,\bar t))}{|\bar x-y^\la|^{n+sp}} dy\leq 0,
 \endaligned \end{equation}
 where $\xi(y,\bar t) $  is between $ \tilde u_\la(\bar x,\bar t)-\tilde u_\la(y,\bar t)$ and $\tilde u (\bar x,\bar t)-\tilde u_\la(y,\bar t)$, $\eta(y,\bar t)$ is between $\tilde u_\la(\bar x,\bar t)-\tilde u (y,\bar t)$ and $\tilde u (\bar x,\bar t)-\tilde u ( y,\bar t)$.
Combining \eqref{eq:o100}, \eqref{eq:o102} and  \eqref{eq:o101}, one has
 \be \label{q20}
 (-\lap)^s_p \tilde u_\la(\bar x,\bar t)-(-\lap)^s_p \tilde {u}(\bar x,\bar t)<0.
 \ee
So
\be \label{20}
\frac{\partial \tilde w}{\partial t}(\bar x,\bar t) +(-\lap)^s_p \tilde u_\la(\bar x, \bar t)-(-\lap)^s_p  \tilde {u} (\bar x, \bar t)<0.
\ee
Since $\tilde w$ and $  w_\la$ have the same sign, one has
 $$w_\la(\bar x,\bar t)= \underset{\bar \Omega \times (0,T]}{\min} w_\la(x,t)<0.$$
 Similarly to the proof of \eqref{q20}, we have
  \be \label{w20}
 (-\lap)^s_p   u_\la(\bar x,\bar t)-(-\lap)^s_p   {u}(\bar x,\bar t)<0.
 \ee
On the other hand, from \eqref{1},  \eqref{w20} and the first inequality of \eqref{1eq:j3}, we have
\begin{eqnarray}\label{q1}
&&\frac{\partial \tilde w}{\partial t}(\bar x,\bar t) +(-\lap)^s_p \tilde u_\la(\bar x,\bar t)-(-\lap)^s_p \tilde u(\bar x,\bar t)\nonumber \\
&=&-m e^{-m  t} w_\la(\bar x,\bar t) +e^{-m  t} \frac{ \partial w_\la}{\partial t}(\bar x,\bar t) +e^{-m(p-1) t}[(-\lap)^s_p u_\la (\bar x,\bar t)-(-\lap)^s_p u (\bar x,\bar t)]\nonumber\\
&\geq & e^{-m  t} [-m w_\la(\bar x,\bar t) +\frac{ \partial w_\la}{\partial t}(\bar x,\bar t)+ (-\lap)^s_p u_\la (\bar x,\bar t)-(-\lap)^s_p u (\bar x,\bar t) ]\nonumber\\
&\geq &   (-m+ c(\bar x,\bar t)) \tilde  w_\la(\bar x,\bar t).
 \end{eqnarray}
Due to $c(x,t)$ is bounded from above, selecting $m$ to make
$$ -m+ c(\bar x,\bar t)<0,$$
then we have
$$ \frac{\partial \tilde w}{\partial t}(\bar x,\bar t) +(-\lap)^s_p \tilde u_\la(\bar x,\bar t)-(-\lap)^s_p \tilde u(\bar x,\bar t) >0,$$
 which contradicts   \eqref{20}.  Hence \eqref{1eq:j4}  is correct.

  ${\rm (ii)}$  If $\Omega$ is unbounded, then   \eqref{1eq:jj4} ensures that the negative minimum of $ w_\lambda(x,t)$  must be reached at some point. Then we can come up with a contradiction based on the same discussion as ($i$).

${\rm (iii)}$ Next we establish \eqref{eq:j6} based on \eqref{1eq:j4}.    Assume that there exists $(x_0, t_0 )\in \Omega \times (0,T]$ such that
$$w_\lambda(x_0, t_0 )=0.$$
Thus $(x_0,t_0)$ is a minimum point of $w_\lambda(x ,t )$. Therefore,
$$ \frac{\partial w_\lambda}{ \partial t}( x_0 ,t_0)\leq0.$$
 Similar to \eqref{eq:o100} and \eqref{eq:o101},
 we obtain
 $$ (-\lap)^{s}_p u_\lambda(x_0,t_0)-(-\lap)^{s}_p u (x_0,t_0)=C_{n,sp}\{J_1+J_2\}~ \mbox{and}~J_2=0.$$
It follows  from the first inequality of \eqref{1eq:j3}  that
 \begin{equation}\label{j14}\aligned
 &   0\leq - \frac{\partial w_\lambda}{\partial t} (x_0,t_0)\\
 &\leq
 (-\lap)^{s}_p u_\lambda(x_0,t_0)-(-\lap)^{s}_p u (x_0,t_0)+c(x_0,t_0)w_\lambda(x_0,t_0)\\
&= C_{n,sp}  \{J_1+J_2\}=C_{n,sp} J_1.
\endaligned\end{equation}
Then we have $J_1\geq 0$. Consequently
$$H(u_\la(  x_0, t_0)-u_\la(y, t_0))-H( {u}(x_0,  t_0)- {u}(y,   t_0))\geq 0$$
and due to the monotonicity of $H$, we derive, for almost all $y\in \Sigma_\la$
$$ [u_\la(x_0,  t_0)-u_\la(y,  t_0)]-[ {u}(x_0,  t_0)- {u}(y,  t_0)]= -w_\la(y,  t_0)\geq 0. $$
Therefore, we have $w_\lambda(x,t_0)\equiv0~a.e.$ in $\Sigma_\lambda,~t_0\in (0,T]$.

Recalling $w_\lambda(x^\lambda,t)=-w_\lambda(x,t),$  we obtain
$$w_\lambda(x,t_0 )= 0,~\mbox{a.e. ~for}~x\in \mathbb R^n,~ t_0\in (0,T].$$

  This completes the proof of Theorem \ref{thmn12}.

\end{proof}

\section{Monotonicity of  the solution for elliptic fractional  $p$-equations }
We first develop a parabolic inequality, then the symmetric and monotonicity for   solutions to   elliptic fractional $p$-equations are derived by the method of moving planes and   maximum principles for the antisymmetric function  to the parabolic fractional $p$-equation. We are equivalent to treating $w_\la (x)$ as $w ( \la,x)$, where the variable $\la$ acts similarly to the variable $t$ in a parabolic equation.

{\bf Proof of Theorem \ref{eqn22}.}  Due to the boundedness of $\Omega$ and convexity with respect to $x_1$  directions, one can assume
$$ \Omega \subset \{ |x_1| \leq b\},~b>0, ~ \partial \Omega \cap \{x_1=-b\} \neq \emptyset$$
and
$$ \Omega_\lambda=\{(x_1,x',\lambda)\in \Omega \mid -b<x_1<\lambda\},~\Sigma_\la = \{ (x_1,x')\in \mathbb R^n \mid x_1<\la \}.$$
 Denote $ u_\la(x)=u(2\lambda -x_1,x')=u(  x^\lambda_1,x')$ and $ w_\la(x)=u_\la(x)-u(x)$.

From equation \eqref{eq:n20191}, since \eqref{eq:n288},    $ (\frac{\partial u_\la}{\partial x_1})(x)=-(\frac{\partial u}{\partial x_1})(x^\la)  $ and $(\frac{\partial u_\la}{\partial x_i})(x)=(\frac{\partial u}{\partial x_i})(x^\la),~i=2,\cdots,n$,  we have
\begin{eqnarray}  \label{30900w}
&&(-\lap)^s_p u_\la(x) -(-\lap)^s_p u (x)   \nonumber \\
&=&h(x^\lambda, u_\la(x), (\nabla  u)(x^\la))-h(x,u(x), \nabla u(x))\\
&=& h(x^\lambda, u_\la , -(\frac{\partial u_\la}{\partial x_1})(x ),(\frac{\partial u_\la}{\partial x_2})(x ),\cdots,(\frac{\partial u_\la}{\partial x_n})(x ) )-h(x,u , \nabla u(x)).\nonumber
\end{eqnarray}
By  \eqref{eq:n288}, if $q_1<0$,  then since Lipschitz continuity, there exists  a positive constant $C_0$ such that
for $ \tau_1>x_1,~x_1+\tau_1<0,$  we have
$$\aligned
&h(\tau_1,x',u,-q_1,q_2,\cdots,q_n)-h(x_1,x',u,q_1,q_2,\cdots,q_n)\\
=& h(\tau_1,x',u,-q_1,q_2,\cdots,q_n)- h(\tau_1,x',u,0,q_2,\cdots,q_n)+h(\tau_1,x',u,0,q_2,\cdots,q_n)\\
&-h(x_1,x',u,0,q_2,\cdots,q_n)
+h(x_1,x',u,0,q_2,\cdots,q_n)-h(x_1,x',u,q_1,q_2,\cdots,q_n)\\
\geq& h(\tau_1,x',u,0,q_2,\cdots,q_n)-h(x_1,x',u,0,q_2,\cdots,q_n)+C_0q_1\\
\geq &C_0q_1.
\endaligned $$
It can be concluded that there is an $L^\infty$ function $\gamma\geq 0$, such that for $\tau_1>x_1,~x_1+\tau_1<0$ and  all $q_1$,
$$ \label{eq:n2000}
h(\tau_1, x',u,-q_1,q_2,\cdots,q_n)-h(x_1,x',u,q_1,q_2,\cdots,q_n)\geq \gamma q_1,
$$
where $\gamma$ depends on $x_1,~\tau_1,~\textbf{q}$.
Therefore, by \eqref{30900w},  for $x\in \Omega_\la,~\la\leq 0$, we derive
$$(-\lap)^s_p u_\la(x) -(-\lap)^s_p u (x)  \geq  h(x,  u_\la(x), (\nabla u_\la)(x)) +\gamma  \frac{\partial u_\la}{\partial x_1} -h(x,u(x), \nabla u(x)), $$
 where $0\leq \gamma \in L^\infty $.

 Due to the Lipschitz continuity of $h(x,u, \textbf{q})$    with respect to  $(u, \textbf{q})$,  for suitable bounded
 functions $a_j(x),~c_\la(x)$, one has
 $$ (-\lap)^s_p u_\la(x)-(-\lap)^s_p u (x)+ \sum^n_1 a_j(x) (w_\la)_j(x) +c_\la(x) w_\la (x)-\gamma \frac{\partial u_\la}{\partial x_1}(x)\geq 0,$$
 here $ c_\la(x)=\frac{h(x,u,\nabla u )-h(x,u_\la,\nabla u )}{u_\la(x)-u(x)},~(w_\la)_j =\frac{\partial w_\la}{\partial x_j}$.
 However $$ \frac{\partial w_\la}{\partial \lambda}(x)=2   u_1 ( x^\lambda_1, x')=-2\frac{\partial u_\la}{\partial x_1}(x),$$
 where $ u_1 ( x^\lambda _1, x')=\frac{\partial u}{\partial x^\la_1}  ( x^\lambda _1, x')$.
Therefore we derive the following parabolic inequality for $w_\la$ as a function of $x$ and $\lambda$, $x\in \Omega_\la$
\begin{equation}\label{2020}
\frac{\gamma}{2} \frac{\partial w_\la}{\partial \lambda}+(-\lap)^s_p u_\la (x)-(-\lap)^s_p u (x) +a_j(x) (w_\la)_j(x) +c_\la (x)w_\la(x)\geq 0.
\end{equation}
In the region $E$ in $(x,\lambda)$ space, we have
$$ E=\{(x_1,x',\lambda)\mid -b <x_1<\lambda <\tilde \lambda, (x_1,x')\in \Omega_\lambda \},~\tilde \la \leq b.$$

Next, we will   apply the method of moving planes to prove  \eqref{eq:n299}.  The proof   consists of the following  two steps.

Step 1.  \textit{Start moving the plane $T_\lambda$  to the right along the $x_1$  direction   from $-b$.}

We will prove that there exists a $ \delta >0$ small enough such that
\begin{equation}\label{eq:mb16}
w_\la(x)\geq0, ~\forall~x\in \Omega_\lambda,~\lambda \in [-b,-b+\delta].
\end{equation}
If \eqref{eq:mb16}  is not valid, denote
$$ B= \underset{-b\leq \la \leq -b+\delta}{\underset{x\in \bar \Omega_\la}{\inf}} w_\la(x)<0.$$
Obviously, $x_1= \la,~ w_\la(x)\equiv0$.
By $u(x)>0, x\in  \Omega$ and $u(x) \equiv 0,~x\in \Omega^c$,   $B$ may obtained for some point $ (\bar \la,\bar x) \in \{(\la,x)\mid (\la,x) \in [-b,-b+\delta]\times \bar \Omega_\la\}$. So one has
$
w_{\bar \la}(\bar x)=B.
$
Noticing that $w_{\bar \la}(x) \geq 0,~x\in \partial \Omega \cap \Sigma_{\bar \la}$, we arrive at $\bar x \in \Omega_{\bar \la}$.
  By the fact that $\Sigma_{-b} \cap \Omega= \partial \Sigma_{-b}\cap \Omega$ for $\la=-b$, one has $w_{-b}(x)\geq 0,~x\in \Sigma_{-b}\cap \Omega$. This implies $\bar \la>-b$. Due to $(\bar \la, \bar t)$ is a minimizing point, we obtain
\begin{equation}\label{23w}
 \nabla_x w_{\bar \la}(\bar x)=0
 \end{equation}
 and
 \begin{equation}\label{22w}
   \frac{\partial w_{\la}(\bar x)}{\partial \la} \mid _{\la=\bar \la} \leq 0.
 \end{equation}
 By calculation, it can be concluded that $ (\partial_{x_1}u)(\bar x^{\bar \la})\leq 0.$
 We have $$ (\nabla_x u_{\bar \la})(\bar x)=(\nabla _x u)(\bar x). $$
 Thus $(w_{\bar \la})_j(\bar x)=0.$

One has  $ w_{- b}(\bar x)\geq 0$ and $w_\la(x)\geq 0,~x\in \Sigma_\la \backslash \Omega_\la,~\la \in [-b,-b+\delta]$.
By \eqref{2020},
we have
 \be \label{000001}
  \frac{\gamma}{2} \frac{\partial w_{\la}( \bar x )}{\partial \lambda}\mid _{\la=\bar \la}+(-\lap)^s_p u_{\bar \la } (\bar x)-(-\lap)^s_p u (\bar x)  +c _{\bar \la}(\bar x)w_{\bar \la}(\bar x)\geq 0.
  \ee
Case 1) When $\gamma=0$, for  $\delta$ small enough, $\Omega_{\bar\la}$ is a narrow region, by the proof of   narrow region principle in \cite{WuN}, it can be seen that
$$ (-\lap)^s_p u_{\bar \la } (\bar x)-(-\lap)^s_p u (\bar x)  +c _{\bar \la}(\bar x)w_{\bar \la}(\bar x)<0.$$
This   contradicts  \eqref{000001} and proves \eqref{eq:mb16}.

Case 2)  When $\gamma>0$, similar to the proof of Theorem \ref{thmn12}, let $\tilde w(x)=e^{-m \la}
w_{\bar \la}(x),~\tilde u_{\bar \la}(x )= e^{-m{  \la}} u_{\bar \la}(x ),~ \tilde u (x )= e^{-m{ \la}} u (x )$, we obtain
$$  \frac{\gamma}{2} \frac{\partial \tilde w ( \bar x )}{\partial \lambda}\mid _{\la=\bar \la}+(-\lap)^s_p \tilde u_{\bar \la } (\bar x)-(-\lap)^s_p\tilde u (\bar x)   \geq (-\frac{\gamma}{2}m-c_{\bar\la}(\bar x))\tilde w(\bar x)>0,$$
 by selecting   $m$ to make   $-\frac{\gamma m}{2}-c_{\bar\la}(\bar x)<0$ due to $c_{\bar\la}(\bar x) $ is bounded.
 On the other hand, by   $\gamma> 0$  and \eqref{20} in Theorem \ref{thmn12}, we have
$$ \frac{\gamma}{2} \frac{\partial \tilde w ( \bar x )}{\partial \la}\mid _{\la=\bar \la}+(-\lap)^s _p\tilde u_{\bar \la} ( \bar x )-(-\lap)^s_p \tilde u  ( \bar x ) < 0.$$
This is a  contradiction.   Hence \eqref{eq:mb16} is right.

Step 2. \textit{Keep moving  the plane to the right until  it reaches its the limiting position $T_{\lambda_0}$ as long as
\eqref{eq:mb16} holds}.
Define $$ \la_0=\sup\{ -b \leq \la \leq 0 \mid  w_\mu(x) \geq 0,~\forall ~x\in \Omega_\mu,~\forall~\mu \leq \la\}.$$
From the continuity of $u(x)$ and the definition of $\la_0$, we derive $$ w_{\la_0}(x)\geq 0,~ \mbox{for ~all}~x\in \Omega_{\la_0}.$$

We conclude that   $\la_0 = 0$ by contradiction arguments.

If $\la_0 = 0$  does not hold, then    $\la_0<0$, we first prove
\begin{equation}\label{800}
w_{\la_0}(x) >0,~x\in \Omega_{\la_0}.
\end{equation}
Otherwise, there exists a point $x_0\in \Omega_{\la_0}$ such that $w_{\la_0}(x_0)= 0.$  So $(\la_0,x_0) $ is a
minimizing point. Since we have  $ (\frac{\partial u}{\partial{x_i}} )(x^{\la_0}_0)=(\frac{\partial u_{\la_0}}{ \partial {x_i} })(x_0)=( \frac{\partial u}{\partial{x_i}})(x_0)$, for $i=2,\cdots,n$ and  $ (\frac{\partial u}{\partial {x_1}})(x^{\la_0}_0)=-(\frac{\partial u_{\la_0}}{\partial{x_1}})(x_0)=-(\frac{\partial u}{\partial  {x_1} })(x_0) \leq 0$ by \eqref{23w}  and \eqref{22w}.  Then we apply property \eqref{eq:n288} of  $h(x,u, \textbf{q})$ and arrive at
$$\aligned
&(-\lap)^s_p u_{\la_0} (x_0)-(-\lap)^s _pu (x_0)\\
=& h(x_0^{\la_0},u_{\la_0}(x_0),- (\frac{\partial u_{\la_0}}{\partial x_1})(x_0), (\frac{\partial u_{\la_0}}{\partial x_2})(x_0), \cdots,   (\frac{\partial u_{\la_0}}{\partial x_n})(x_0))- h(x_0 ,u (x_0), \nabla u (x_0))\\
\geq &   h(x_0 ,u_{\la_0} (x_0), (\frac{\partial u }{\partial x_1})(x_0),\cdots,   (\frac{\partial u }{\partial x_n})(x_0)  )- h(x_0 ,u (x_0), \nabla u (x_0))
=0,\endaligned$$
which contradicts
$$(-\lap)^s_p u_{\la_0} (x_0)-(-\lap)^s _pu (x_0) =C_{n,sp}\{I_1+I_2\} <0 $$
  by \eqref{eq:o100} with $I_1<0,~I_2=0$ and  $w_{\la_0} \not \equiv 0$.

Suppose that $\la_0<0$, it can be inferred that there exists $ \epsilon  _0 > 0$ small enough such that
\begin{equation}\label{389}
w_\la(x)\geq 0,~x\in \Omega_\la,~\forall ~\la \leq \la_0+\epsilon_0.
\end{equation}
Suppose \eqref{389} is incorrect, we have
$$ B_k=\underset{\la_0\leq \la \leq \la_k}{\underset{x\in \bar \Omega_\la}{\inf}} w_\la(x)<0,~\mbox{for~a~sequence~of ~} \la_k \searrow \la_0,~\mbox{as}~k \rightarrow +\infty.$$
The minimum $B_k$ can be obtained for some $\mu_k \in (\la_0, \la_k], x_k \in \Omega_{\mu_k}
 $ where
$w_{\mu_k} (x_k) = B_k$ by the same reason as in Step 1. So
\be \label{26}
   ( w _{\mu_k})_j(x_k)=0.
  \ee
   By the property of $\la_k$, we have $\mu_k \rightarrow \la_0$ as $k\rightarrow +\infty$.

\textbf{Case 1)}  If  $ \gamma>0$,
from \eqref{2020} and \eqref{26}, let $\tilde w(x)=e^{-m\la} w_{\mu_k}(x),~m>0$, we obtain
\begin{equation}\label{1223}
 \frac{\gamma}{2} \frac{\partial \tilde w }{\partial \la}\mid_{\la=\mu_k} +(-\lap)^s_p \tilde u_{\mu_k} ( x_k)-(-\lap)^s_p \tilde u  ( x_k) \geq (-\frac{\gamma}{2}m-c_{\mu_k} ( x_k)) \tilde w (  x_k)>0,
 \end{equation}
 by selecting $m$ to make $ -\frac{\gamma}{2}m-c_\la ( x_k)<0$ because of the boundedness of  $c_{\mu_k}  (x)$ and $\gamma$.
On the other side, from \eqref{22w}, $\gamma> 0$  and \eqref{20} in Theorem \ref{thmn12}, one has
$$ \frac{\gamma}{2} \frac{\partial \tilde w }{\partial \mu_k}( x_k)+(-\lap)^s _p\tilde u_{\mu_k} ( x_k)-(-\lap)^s_p \tilde u  ( x_k) < 0.$$
This contradicts \eqref{1223},  thus we   have \eqref{389}.  This is a contradiction with the definition of $\la_0$.

\textbf{Case 2)} If $\gamma=0$, from  \eqref{2020} and \eqref{26}, we have
$$  (-\lap)^s_p   u_{\mu_k}  ( x_k)- (-\lap)^s_p   u ( x_k)+ c_{\mu_k}  ( x_k)    w _{\mu_k}(  x_k) \geq 0.$$
It follows from  \eqref{800} that for any $\delta>0$
$$w_{\la_0}(x)>c_0>0,~\forall~x\in \Omega_{\la_0-\delta}.$$
By the continuity of $w_\la$ with respect to $\la$,  for $ \la_0<\mu_k \leq \la_k$, such that
\be \label{27}
w_{\mu_k}(x) \geq 0,~\forall~x\in \Omega_{\la_0-\delta},~\forall~\mu_k \in (\la_0, \la_k).
\ee
For narrow region $\Omega_{\mu_k} \backslash \Omega_{\la_0-\delta}$,   applying    narrow region principle in \cite{WuN},   we have
$$ w_{\mu_k}(x) \geq 0,~\forall~x\in \Omega_{\mu_k} \backslash \Omega_{\la_0-\delta}.$$
This together with \eqref{27} implies
$$ w_{\mu_k}(x) \geq 0,~\forall~x\in \Omega_{\mu_k} ,~\forall ~\mu_k \in (\la_0,\la_0+\epsilon_0) .$$
 Therefore this contradicts the definition of   $\la_0$.

Combining \textbf{Case 1)} and \textbf{Case 2)}, we   have
$\la_0=0$.
Hence
\begin{equation}\label{00}
u(x_1,x') \leq u(-x_1,x'),~\forall~(x_1,x')\in \Omega,~x_1<0.
\end{equation}
Furthermore, similar to the proof of \eqref{800}, we can actually deduce that
$$ w_\la (x)>0,~x\in \Omega_\la,~\forall~\la<0.$$
For any $(x_1,x'),(\tilde x_1, x')\in \Omega$ with $0>x_1>\tilde x_1$, we take $\la=\frac{x_1+\tilde x_1}{2}$. Then we obtain
$$u(x_1,x')>u(\tilde x_1,x')$$
and thus $u(x_1, x')$ is strictly increasing in the left half of $\Omega$
 in $x_1$-direction.

 In addition ,  if $h(x_1,x',u,q_1,q_2,\cdots,q_n)=h(-x_1,x',u,-q_1,q_2,\cdots,q_n)$,
 then we derive $\bar u (x_1,x')=u(-x_1,x')$ also solves \eqref{eq:n20191}. Thus we have obtained that
 $$ \bar u(x_1,x')\leq \bar u(-x_1,x'),~\forall~(x_1,x')\in \Omega,~x_1 <0,$$
or equivalently,
$$   u(x_1,x')\geq   u(-x_1,x'),~\forall~(x_1,x')\in \Omega,~x_1 <0.$$
Combining the above formula  with \eqref{00},  it can be concluded that
$$   u(x_1,x')=   u(-x_1,x'),~\forall~(x_1,x')\in \Omega,~x_1 <0,$$
$i.e.$, $u$ is symmetric in the $x_1$ direction about $ \{x_1 = 0\}$.
  Theorem  \ref{eqn22} has been proven.

\section{ Monotonicity and symmetry of the solution for parabolic fractional $p$-equations }
In this section, we use the  maximum principle for antisymmetric functions to parabolic  fractional $p$-equations to  derive monotonicity and symmetry of fractional parabolic equations in $\mathbb R^n$ and the bounded domain.
\subsection{The bounded domain }
{ \bf Proof of Theorem \ref{thmn1}.} By \eqref{eq:n2}, we derive that  the initial values of $u$ are not a constant and  strictly increasing in $x_1$ when $x_1<0$, then
as time $t$ goes on, we only need  to prove that $u$    is monotonicity increasing in $x_1$ for each time when $x_1<0$.
Set $$-b:=\min\{x_1 \mid x\in \Omega\},~b>0,  ~\Omega_\lambda=\{x\in \Omega\mid -b<x_1<\lambda \}$$
and $$ w_\lambda (x,t)=u(x^\lambda,t)-u(x,t).$$
We will show that
\begin{equation}\label{eq:n7}
w_\lambda(x,t)\geq 0,~x\in \Omega_\lambda,~t\in (0,T),~\lambda\in (-b,0).
\end{equation}
Since \eqref{20191} and $u(x^\lambda,t)$ is also a solution of \eqref{eq:n1}, we obtain
\begin{eqnarray}\label{eq:n8}
&&\frac{\partial w_\lambda}{\partial t}(x,t)+(-\lap)^s _pu_\lambda (x,t)-(-\lap)^s _pu (x,t)\nonumber\\
&=&F(x^\lambda,t,u(x^\lambda,t))-F(x,t,u(x,t))\nonumber\\
&\geq& F(x ,t ,u(x^\lambda,t ))-F(x,t ,u(x ,t ))\nonumber\\
&:=  &c_\la(x ,t ) w_\lambda(x ,t ),~x\in \Omega_\lambda,~0<t<T,
 \end{eqnarray}
where $ c_\la(x,t )=\frac{F(x ,t ,u(x^\lambda,t ))-F(x ,t ,u(x ,t ))}{u(x^\lambda,t )-u(x ,t )} $ is bounded by
 Lipschitz continuity of $F$ with respect to  $u$.

Due to the initial  assumption \eqref{eq:n2}, we arrive at
$$w_\lambda(x,0)\geq 0,~x\in \Omega_\lambda, ~\lambda \in (-b,0).$$
For each fixed $\lambda$,   we show that   \eqref{eq:n7}.
From exterior condition \eqref{eq:n3}, we obtain
$$ w_\la(x,t)\geq 0,~(x,t) \in (\Sigma_\la\backslash \Omega_\la)\times [0,T).$$
Applying  Theorem  \ref{thmn12}, this verifies \eqref{eq:n7}.
 So we derive
 $$ w_\lambda(x,t)\geq 0,~x\in \Omega_\lambda,~\lambda\in (-b,0),~0 \leq t \leq T,$$ which is equivalent to
\begin{equation}\label{2020w}
u(x_1,x',t)\leq u(-x_1,x',t)~\mbox{ for}~
(x_1,x')\in \Omega,~x_1<0~\mbox{and }~0\leq t \leq T.
\end{equation}

Next we prove  \eqref{2q2q}.
By condition  \eqref{w20191}, $\hat u(x_1,x',t)=u(-x_1,x',t)$ is also the solution of \eqref{eq:n1}.
Since the initial-exterior conditions of $u$  are symmetric about the plane $\{x_1=0\}$, then
 apply the same process to the function
$\hat u(x_1,x',t)$, one has
$$\hat u(x_1,x',t) \leq \hat u(-x_1,x',t),~x_1<0,~t\in [0,T],$$
$i.e.,$
\begin{equation} \label{eq:o331}
  u(x_1,x',t)\geq   u(-x_1,x',t),~x_1<0,~t\in [0,T].
\end{equation}
 Combining \eqref{2020w} and \eqref{eq:o331},   we have
 $$ u(x_1,x',t)=u(-x_1,x',t) ~\mbox{for} ~ x_1\leq 0,~t\in [0,T].$$
 Hence we obtain
 \begin{equation}\label{eq:o34}
 \frac{\partial w_\lambda}{\partial x_1}(x,t) \leq 0,~x\in \Omega_\lambda,~t\in (0,T],~\lambda\in (-b,0].
 \end{equation}
 From the definition of $w_\lambda(x,t)$, it yields
 $$\frac{\partial u}{\partial x_1}(x,t)\geq 0,~x\in \Omega,~x_1<0,~0<t\leq T.$$
Therefore,
 $$ u(x_1,x',t)=u(-x_1,x',t) ~\mbox{and}~\frac{\partial u}{\partial x_1}(x,t)\geq 0~\mbox{for}~x\in \Omega, x_1<0,~ t\in (0,T].$$

Now we show that
$$\frac{\partial u}{\partial x_1}(x,t)>0,~x_1<0,~t\in (0,T].$$
We have proven that
$$w_\lambda(x,t)\geq 0,~x\in \Omega_\lambda,~0< t\leq T,~ \lambda \leq 0 $$
and \begin{equation}\label{eq:n201914}
w_0(x,t)= 0,~x\in \Omega_0,~t\in (0,T] .
\end{equation}
For fixed $\lambda<0$,   we  conclude that
 \begin{equation}\label{eq:n201915}
 w_\lambda(x,t)> 0,~x\in \Omega_\lambda,~t\in (0,T], ~\lambda \in[-b, 0).
 \end{equation}

Denote $$\tilde u(x,t)=e^{-mt}u(x,t),~ \tilde w_\lambda(x,t)=e^{-mt}w_\lambda(x,t),~m>0.$$
From $\tilde w_\la$ and $w_\la$ have the same sign,   we conclude that \eqref{eq:n201915} by proving
 \begin{equation}\label{1eq:n201915}
 \tilde w_\lambda(x,t)> 0,~x\in \Omega_\lambda,~0< t\leq T, ~\lambda \in[-b, 0).
 \end{equation}
If \eqref{1eq:n201915} is wrong,  there exists $\bar x\in \Omega_\lambda$ and the first $\bar t\in (0,T]$ such that
$$\tilde w_\lambda (\bar x, \bar t)=\underset{\bar \Omega_\lambda \times (0,T]}{\min}\tilde w_\lambda (x,t)=0.$$
 From condition \eqref{eq:n2} and \eqref{eq:n3}, it can be further inferred that $(\bar x,\bar t)$ is in the interior of $\Omega_\lambda\times (0,T]$. On the one hand,  by \eqref{eq:n8}, setting $m=2|c(x,t)|$, we have
\begin{equation}\label{eq:n201913}
\frac{\partial \tilde w_\lambda}{\partial t}(\bar x,\bar t)+(-\lap)^s_p \tilde u_\lambda (\bar x,\bar t)-(-\lap)^s_p \tilde u (\bar x,\bar t) \geq (-m+c(\bar x,\bar t))\tilde w_\lambda(\bar x,\bar t)\geq0.
\end{equation}
On the other hand, by $\frac{\partial \tilde w_\lambda}{ \partial t}(\bar x,\bar t)\leq 0$,
similar to \eqref{20} in  the proof of  Theorem  \ref{thmn12}, we derive
$$\aligned
\frac{\partial \tilde w_\lambda}{ \partial t}(\bar x,\bar t)+(-\lap)^s_p\tilde u_\lambda(\bar x,\bar t)-(-\lap)^s_p\tilde u(\bar x,\bar t)<0,
\endaligned $$
which contradicts \eqref{eq:n201913}.
So \eqref{eq:n201915} is correct. Combining \eqref{eq:n201914} and  \eqref{eq:n201915}, we have
$$\frac{\partial w_\lambda}{\partial x_1}(x,t)<0,~x_1<0,~0\leq t\leq T.$$
From the definition of $w_\lambda$, we obtain
$$\frac{\partial u}{\partial x_1}(x,t)>0,~x_1<0,~0\leq t\leq T.$$
 Therefore \eqref{eq:n201912} is  right.

 {\bf Proof of Theorem \ref{thmn3}.}    Choosing any   direction as the $x_1$ direction. Based on the monotonically decreasing property of  $F(x,t,u(x,t))=F(|x|,t,u(x,t))$ with respect to $|x|$   and the Lipschitz continuity of $F$  with respect to $u$ for $t\geq t_0$, due to \eqref{eq:nll13} and \eqref{eq:ll112},
 similar to the proof of Theorem \ref{thmn1},     we arrive at
 $$ u(x_1,x',t)=u(-x_1,x',t),~x_1<0,~ t\in (t_0 , T].$$
 Namely, $$w_0(x,t)\equiv 0,~x\in B_1(0),x_1\leq 0,~t\in (t_0 , T].$$
 Due to the direction of $x_1$ can be selected arbitrarily, we have proven that
 \eqref{eq:n201918}.  The proof of monotonicity is similar to the proof of  \eqref{eq:n201915}, and we will omit it here. This completes the proof.

\subsection{The whole space }
In this section, we will prove Theorem \ref{thmn4}.

 {\bf Proof of Theorem \ref{thmn4}.}  Choosing any direction as the $x_1$ direction.
   Set $$ w_\lambda (x,t)=u(x^\lambda,t)-u(x,t)$$
 and
$$ \Sigma_\lambda =\{ x\in \mathbb R^n \mid x_1 <\lambda <0\}.$$
\textit{Step 1. }   We will show that for sufficiently negative $\la$, it holds
\be \label{4281}
w_\la(x,t)\geq 0,~x \in \Sigma_\la,~  t\in(0,T].
\ee
Due to $u(x^\lambda,t)$ is also a solution to \eqref{eq:wn1}, similar to \eqref{eq:n8},  by \eqref{220191}  we have
\begin{equation}\label{11eq:n8}\aligned
 \frac{ \partial w_\lambda}{\partial t}(x,t)+(-\lap)^s_p u_\lambda (x,t)-(-\lap)^s_p u  (x,t)
 \geq   c_\la(x ,t ) w_\lambda(x ,t ),~x\in \Sigma _\lambda,~t\in (0,T],
\endaligned\end{equation}
where $ c_\la(x,t )=\frac{F(x ,t ,u(x^\lambda,t ))-F(x ,t ,u(x ,t ))}{u(x^\lambda,t )-u(x ,t )} $ is bounded by  Lipschitz continuous of $F$  with respect to   $u$.

Now we analyze the initial   condition and the asymptotic behavior of $u$,
By \eqref{xeq:n2}, it knows that  the initial values of $u$ is strictly increasing in $x_1$ when $x_1<0$ and is not a constant. Thus we obtain
 \begin{equation}\label{x567}
 w_\lambda(x,0)\geq 0,~x\in \Sigma_\lambda,~\lambda \leq 0 .
 \end{equation}
 For  fixed $\lambda$ and all $t\in[0,T]$, the condition \eqref{4283}  hints that
$$u(x,t)\rightarrow0,~\mbox{as}~|x|\rightarrow+\infty.$$
Because $|x^\lambda|\rightarrow+\infty$, as $|x|\rightarrow+\infty$,  we derive
$$u_\lambda(x,t)=u(x^\lambda,t)\rightarrow0,~0\leq t\leq T. $$
Thus we obtain
\begin{equation}\label{428300}
w_\lambda(x,t)\rightarrow 0,~\mbox{as}~|x|\rightarrow+\infty,~0\leq t\leq T.
\end{equation}
From \eqref{x567} and \eqref{428300}, it can be inferred that the point where $w_\la<0$ will not be reached at $t=0$ and $| x |\rightarrow +\infty$. Due to the unknown sign of term $c_\la(x)$ in inequality \eqref{11eq:n8}, we adopt the proof method of the maximum principles of antisymmetric functions (Theorem \ref{thmn12}).
Set
 $$\tilde u_\lambda(x,t)=e^{-mt}u_\lambda (x,t),\tilde u(x,t)=e^{-mt}u(x,t), \tilde w_\lambda(x,t)=e^{-mt}w_\lambda (x,t),~m>0.$$
To prove \eqref{4281}, we only need claim that
\be \label{4282}
\tilde w_\la(x,t)\geq 0,~x \in \Sigma_\la,~  t\in(0,T].
\ee
Otherwise,  by \eqref{x567}, \eqref{428300}  and $w_\la(x,t)\equiv0,~x\in T_\la$, there   exists $(x^0,t^0)\in \Sigma_\la \times (0,T]$such that
$$ \tilde w_\la (x^0,t^0)= \underset{\Sigma_\la \times (0,T]}{\min}\tilde w_\la(x,t)<0.$$
 By \eqref{11eq:n8},  we  obtain
 $$  \frac{\partial \tilde w_\lambda}{\partial t}(x^0,t^0)+(-\lap)^s_p \tilde u_\lambda (x^0,t^0) -(-\lap)^s_p \tilde u  (x^0,t^0)\geq (-m+c_\la(x^0,t^0)) \tilde w_\la (x^0,t^0).$$
Because of the boundedness of $c_\la(x,t)$, we select $ m$ to make $ c_\la(x^0,t^0)-m<0$. So
 \begin{equation}\label{20201w}
 \frac{\partial \tilde w_\lambda}{\partial t}(x^0,t^0)+(-\lap)^s \tilde u_\lambda (x^0,t^0)-(-\lap)^s \tilde u  (x^0,t^0)>0.
 \end{equation}
On the other side, similar to the proof process of \eqref{20}   in    Theorem  \ref{thmn12}, we derive
$$\frac{\partial \tilde w_\lambda}{\partial t}(x^0,t^0)+(-\lap)^s_p\tilde u_\lambda(x^0,t^0)-(-\lap)^s_p \tilde u  (x^0,t^0) <0,$$
which contradicts \eqref{20201w}. Thus we have \eqref{4282}. Step 1 is then completed.

\textit{Step 2.} Let
$$ \la_0:=\{ \la <0\mid w_\mu(x,t)\geq 0,\forall x\in \Sigma_\mu,t\in (0,T],\forall \mu \leq \la\}.$$
Our   goal is to prove $\la_0=0$. Suppose that $\la_0<0$, then there exists a sufficiently small $\varepsilon>0$ such that $\la_0+\varepsilon<0$. We  will show that $w_\la(x,t) \geq 0,(x,t)\in \Sigma_\la \times (0,T] $ for $\la=\la_0+\varepsilon$, which contradicts the definition of $\la_0$.  If $w_{\la_0+\varepsilon}<0$ somewhere in $\Sigma_{\la_0+\varepsilon}\times (0,T]$,
by the asymptotic behavior of $u$ \eqref{428300} and the initial   condition \eqref{x567}, one gets that the negative minimum point is obtained in the interior of $\Sigma_{\la_0+\varepsilon}\times (0,T]. $ By the same argument as that in Step 1, we know that it is impossible. Therefore, we claim that $\la_0=0$,
 namely,
\begin{equation}\label{x2020w}
u(x_1,x',t)\leq u(-x_1,x',t)~\mbox{for}~x_1\leq 0,~  t \in[0,   T].
\end{equation}

In addition, because $F(x,t,u)=F(|x|,t,u)$, $\hat u(x_1,x',t)=u(-x_1,x',t)$ is also a solution of \eqref{eq:n1}.
Under the condition that the initial value of $u$ are  radially symmetric with respect to  the origin and decreasing in $|x|$, and after the same discussion of the function
$\hat u(x_1,x',t)$, we   obtain
$$\hat u(x_1,x',t) \leq \hat u(-x_1,x',t),~x_1<0,~0\leq t\leq T.$$
That is
\begin{equation} \label{xeq:o331}
  u(x_1,x',t)\geq   u(-x_1,x',t),~x_1<0,~0\leq t\leq T.
\end{equation}
Alternatively, we move the plane from positive infinity to the left to obtain  \eqref{xeq:o331}.
 Combining \eqref{x2020w} and \eqref{xeq:o331},   we have
 $$ u(x_1,x',t)=u(-x_1,x',t) ~\mbox{for} ~ x_1\leq 0,~t\in [0,T].$$

     Due to to the arbitrariness of the $x_1$  direction   can be selected arbitrarily, we conclude that $u(x,t)$ is radially symmetric about the origin for each $0\leq t\leq T$.  The monotonicity of $u(x,t)$  is similar to   the discussion   of   \eqref{eq:n201915}.

 This completes the proof of Theorem \ref{thmn4}.

\vskip 6mm
\noindent{\bf Acknowledgements}

\noindent
The   author was supported by the   National Natural Science Foundation of China (No. 12101530), Scientific and Technological Key Projects of Henan Province (No.232102310321) and Nanhu Scholars Program for Young Scholars of XYNU (No. 2023).

\end{document}